\documentclass{article}
\usepackage{amssymb} 
\usepackage{a4} 
\usepackage{amsmath}
\usepackage{latexsym}

\newcommand{\bb}[1]{{\mathbb #1}}    

\newcommand{\bbR}{{\bb R}}

\newcommand{\bbC}{{\bb C}}

\newcommand{\lmod}{\backslash}

\newcommand{\SU}{{\rm SU}} 
\newcommand{\U}{{\rm U}}

\newcommand{\tr}{{\rm tr}}

{ \par  \medskip  \noindent  {\bf Definition} \hspace{0.5em} }%
{\par \medskip } 

\newtheorem{example1}{Example}[section]

\newenvironment {example}%
{ \begin{example1} \rm }%
{ \end{example1}} 
\newenvironment {proof}%
{ \noindent {\em Proof. }}%
{\hspace*{\fill}$\Box$\par \medskip } 
\newenvironment{prf}[1]%
{ \noindent {\em #1 \hspace{0.5em}} }%
{\hspace*{\fill}$\Box$\par \medskip } 
{{\em Remark \hspace{0.5em}} }%
{\par \medskip }

\newtheorem{proposition}{Proposition}
\newtheorem{definition1}[proposition]{Definition}
\newtheorem{theorem}[proposition]{Theorem}
\newtheorem{lemma}[proposition]{Lemma}
\newtheorem{corollary}[proposition]{Corollary}
\newenvironment{definition}%
{ \begin{definition1} \rm }%
{ \end{definition1}}

\newfam\eufam%
\font\euzw=eufm10 scaled 1200%
\font\euac=eufm9%
\textfont\eufam=\euzw\scriptfont\eufam=\euac%
\def\fr{\fam\eufam\euzw}%
\newcommand{\Sr}{{\rm S}}  
\newcommand{\CP}{\bbC {\rm P}}  
\newcommand{\CH}{\bbC {\rm H}}  
\renewcommand{\Im}{{\rm \, Im}}
\renewcommand{\Re}{{\rm \, Re}}

\title{Proper Affine Hyperspheres which fiber over
Projective Special K\"ahler Manifolds}

\author{Oliver Baues \thanks{e-mail: oliver@math.ethz.ch} \\
Departement Mathematik\\ 
ETH-Zentrum \\ 
R\"amistrasse 101 \\
CH-8092 Z\" urich 
\and
Vicente Cort\' es \thanks{e-mail: vicente@math.uni-bonn.de} \\
Mathematisches Institut \\
Universit\"at Bonn\\
Beringstra\ss e 1 \\
D-53115 Bonn
}

\date{May  17, 2002}

\begin{document}

\maketitle
\begin{abstract} \noindent 
We show that the natural $S^1$-bundle over a 
projective special  K\"ahler manifold 
carries the geometry of a proper affine hypersphere 
endowed with a Sasakian structure. The construction generalizes the 
geometry of the Hopf-fibration $\Sr^{2n+1} \longrightarrow \CP^n$ in
the context of projective special  K\"ahler manifolds. As an application
we have that a natural circle bundle over the Kuranishi moduli
space of a Calabi-Yau threefold is a Lorentzian proper affine
hypersphere. 
\end{abstract} 

\bigskip \bigskip 
\section*{Introduction}
In a previous paper \cite{bc1}, the authors proved 
that any simply connected special K\"ahler manifold 
admits a canonical immersion into affine space 
as a parabolic affine hypersphere. 
A particular important class of special K\"ahler manifolds 
are {\em conic\/} special K\"ahler manifolds. These are
by definition special K\"ahler manifolds which are locally
modelled on a complex cone over some complex projective
manifold which is then called a {\em projective special K\"ahler manifold}. 
The purpose of this paper is to provide an understanding
of the particular (affine) differential geometry which is 
canonically associated with projective special K\"ahler manifolds. 

Whereas the conic special K\"ahler manifold $M$ which
is associated with a simply connected 
projective special K\"ahler manifold $\bar{M}$
carries the geometry of a {\em parabolic\/} (or improper) affine 
hypersphere, we show that the total space $S$ of a natural 
circle bundle $S \rightarrow \bar{M}$ is a {\em proper\/} affine hypersphere.
The $S^1$-action on $S$ induces a
Sasakian structure on $S$ which is compatible 
with the affine differential geometry in a very specific sense. 
Moreover, all information about the conic special K\"ahler geometry on
$M$ is encoded in the {\em affine Sasakian geometry\/} on $S$. 

Lu showed  \cite{Lu} that every {\em complete\/} affine special 
K\"ahler manifold with a positive definite metric is 
flat. Using a well known result of Calabi \cite{Cal2} on
complete convex affine hyperspheres we obtain  an analogous result 
for {\em projective\/} special K\"ahler manifolds: We show 
that if $\bar{M}$ is a (simply connected) complete 
projective special K\"ahler manifold with a definite affine metric
on $S$ then $\bar{M}$ is isometric to $\CP^n$ with
the canonical Fubini-Study metric.  

The construction of the affine sphere $S$ over a projective 
special K\"ahler manifold naturally
relates to well known canonical data on the Kuranishi
moduli space for Calabi-Yau three-manifolds. Thereby 
we show that a natural circle bundle over the 
Kuranishi moduli space admits a canonical structure
of an affine hypersphere with affine metric of Lorentzian  
signature. 

If $\bar{M}$ is complete, and the metric on $S$ is not 
definite, as in the case of  Kuranishi
moduli spaces, then interesting complete models for projective 
special K\"ahler manifolds do exist. We describe all fibrations 
$S \rightarrow \bar{M}$ which admit a transitive   
semisimple group of automorphisms preserving the projective 
special K\"ahler structure on the base $\bar{M}$. These are  
particular examples of homogeneous Lorentzian affine hyperspheres
fibering over Hermitian symmetric spaces.

\section{Preliminaries}
\label{prelimSec}

\subsection{Affine hypersurfaces}
For the convenience of the reader, and to fix the
notation, we recall the basic definitions of affine
differential geometry of hypersurfaces in $\bbR^{n+1}$  
and the definition of {\em affine hyperspheres}. 
For more details, see for example \cite{NS,Cal2}.
Let $\det$ denote the standard volume form on $\bbR^{n+1}$,
and $\nabla$ the standard flat connection on $\bbR^{n+1}$. 
In the context of {\em affine immersions\/} we
consider manifolds with a semi-Riemannian metric $g$ and  
a torsionfree connection $\hat{\nabla}$
so that 
\begin{itemize}
\item[i)] the cubic tensor $\hat{\nabla} g$
is totally symmetric, and 
\item[ii)] the metric
volume form $\theta_g$ is $\hat{\nabla}$-parallel.
\end{itemize} 
The data  $(g,\hat{\nabla})$ are then said to satisfy the
{\em compatibility condition\/} i) and the {\em equiaffine condition\/} ii).
Every nondegenerate hypersurface immersion $\psi: M \rightarrow \bbR^{n+1}$  
induces data $(\hat{\nabla}, g)$ on $M$ which satisfy i) and ii) via
the fundamental formula
\begin{equation} \label{Gformula}
   \nabla_X Y = \hat{\nabla}_X Y + g(X,Y) E \;, 
\end{equation}
where $X,Y$ denote vector fields on $M$, and
$E$ is the {\em affine normal\/} of the immersion. 
(Note that the notation identifies $M$ as a submanifold 
of $\bbR^{n+1}$.)
The affine normal $E$ is a canonical normal vector field along $\psi$
which is defined up to sign by the condition that the pair $(\hat{\nabla},g)$
satisfies ii), and the
normalizing condition 
\begin{itemize}
\item[iii)]
$\det(E, \ldots) = \theta_g$ on $M$.
\end{itemize} 
The metric $g$ is then called the {\em Blaschke metric\/} and the immersion 
$\psi :(M,\hat{\nabla} , g) \rightarrow \bbR^{n+1}$ a {\em Blaschke immersion}. 
The tensor $A= - D E$ is horizontal along $\psi$ 
and is called the {\em shape tensor\/} of the immersion. 
The quantity $H = {1 \over n} \tr A$ is 
called the {\em affine mean curvature}.
If $\hat{\nabla}$ is flat and $n>1$ then, by the equation of Gau\ss, $A=0$ and the 
affine normal is the restriction of a 
constant vector field. In this case, $\psi$
is called a {\em parabolic\/} (or improper) {\em affine hypersphere}. 
If the shape tensor equals a constant multiple of the identity, $A = \kappa \, id$,
where $\kappa \neq 0$, $\psi$ is called a {\em proper affine hypersphere}. 
In this case, $\hat{\nabla}$ is projectively flat. An affine hypersphere has
constant mean curvature $H= \kappa$. 

Let $M$ be a manifold with data $(\hat{\nabla},g)$ which satisfy i) and ii).
We put $\hat{\nabla}^*$ for the conjugate connection of $\hat{\nabla}$ with
respect to $g$. It is torsionfree by the compatibility condition i). 
Then the fundamental theorem of affine differential geometry 
asserts that a simply connected manifold $M$ with data $(\hat{\nabla},g)$
arises from a  Blaschke immersion $\psi$ if and only if 
the {\em integrability condition}
\begin{itemize}
\item[iv)] ${\hat{\nabla}}^*$ is projectively flat
\end{itemize} 
is satisfied. The immersion $\psi$ is determined by the data
$(\hat{\nabla},g)$ up to composition with an unimodular affine 
transformation. A special case arises if $\hat{\nabla}$ is flat.  
Then it is easily seen that $\hat{\nabla}^*$ is also flat. Hence, iv) is
satisfied and $M$ is a parabolic affine hypersphere. We also mention
that the data $(\hat{\nabla},g)$ arise from an immersion as an affine sphere
if and only if the cubic tensor  $C=\hat{\nabla} g$ has totally
symmetric derivative  $\hat{\nabla} C$. If $(M, \hat{\nabla},g)$ 
is a manifold which satisfies the integrability conditions for 
a Blaschke immersion as an affine sphere we say that 
$M$ has the structure of an affine sphere. 

\subsection{Special K\"ahler manifolds} \label{SK}
We recall some basic notions and constructions from special K\"ahler
geometry. For more details the reader can consult \cite{ACD}, 
and also \cite{Freed}.
A {\em  special K\"ahler manifold} $(M,J,g,\nabla )$ is a  
(pseudo-) K\"ahler manifold $(M,J,g)$ together with a flat torsionfree 
connection $\nabla$ such that $\nabla \omega = 0$, where 
$\omega = g(\cdot , J \cdot )$ is the  K\"ahler form, and 
such that $\nabla J$ is symmetric, i.e.\ $d^\nabla J(X,Y) := (\nabla_XJ)Y - 
(\nabla_YJ)X = 0$ for all vector fields $X$ and $Y$.
 
More precisely, one should speak of {\em affine\/}  special K\"ahler manifolds 
since there is also the notion of a projective special K\"ahler manifold. 
In fact, there is a class of (affine)  special K\"ahler manifolds $(M,J,g,
\nabla)$, 
which are called {\it conic special K\"ahler manifolds} and 
which are characterized by the existence of a  
local holomorphic\linebreak[4] $\bbC^*$-action 
$\varphi_\lambda : M \rightarrow M$, 
$\lambda = re^{it}\in \bbC^*$, with the property: 
\[ (\varphi_\lambda)_* X = 
r\cos t X + r\sin tJX \]  
for all $\nabla$-parallel vector fields $X$ on $M$. 
Under appropriate regularity assumptions on the action, the projection 
$$ \pi:M \longrightarrow \bar{M}=P(M) $$ 
onto the space of orbits $\bar{M}=P(M)$ is a
holomorphic submersion onto a complex (Hausdorff-) manifold. Then
$\bar{M}$ inherits a \mbox{(pseudo-)} K\"ahler metric $\bar{g}$ from $(M,g)$, 
and the base $(\bar{M},\bar{g})$ is called a 
{\it projective special K\"ahler manifold}. Although, 
strictly speaking, the fully fledged projective special K\"ahler geometry 
is encoded in the geometric data on the bundle $\pi:M \rightarrow \bar{M}$. 

Special K\"ahler manifolds may also be characterized in
terms of {\em complex Lagrangian immersions\/} (see \cite{ACD}). 
In fact, any simply connected special K\"ahler manifold $(M,J,g,\nabla )$
has a canonical realization as a (pseudo-) K\"ahlerian 
immersed Lagrangian submanifold of a pseudo-Hermitian, 
complex symplectic vector space $(V,\gamma, \Omega )$ with 
split signature. 
This means that there exists a holomorphic Lagrangian 
immersion $\lambda: M \rightarrow V$ so that 
$g = \lambda^* \gamma$ is the pull-back of the hermitian 
product $\gamma$.  
Moreover, 
the projection  
onto the subspace $V^{\tau}$  of real points for the real structure 
$\tau$ defined by the relation 
$\Omega= -i \, \gamma(\cdot, \tau \cdot)$   
gives local flat coordinates on $M$ which determine 
the flat connection $\nabla$.  
The holomorphic 
Lagrangian immersion $\lambda$ is determined 
by the data $(g,\nabla )$ up to a complex 
affine transformation which preserves $\gamma$ and $\Omega$. 
Conic special K\"ahler manifolds may be realized 
by immersions $\lambda$ which are equivariant with 
respect to the natural $\bbC^*$-action on $V$. 
$\lambda$ is then uniquely
determined up to a complex 
linear transformation which preserves $\gamma$ and $\Omega$. 
We then call $\lambda$ a {\em compatible\/} Lagrange immersion
of the (conic) special  K\"ahler manifold $M$.

\section{The local geometry} \label{localgeometry}
It is well known that holomorphic Lagrangian immersions 
$\lambda$ into a complex  $2n$-dimensional 
symplectic vector space $V$ 
are locally of the form $\lambda = \lambda_F := 
dF: U \rightarrow T^*\bbC^n \cong V$, 
where $F$ is a holomorphic function defined on some domain 
$U \subset \bbC^n$.  
The K\"ahler condition for the holomorphic Lagrangian 
immersion $\lambda_F$ is an open condition on the real 2-jet
of $F$. Conic special K\"ahler manifolds correspond to 
potentials which are homogeneous of degree 2. 
Therefore the local geometry of (conic)     
special K\"ahler manifolds 
may be described in terms of a holomorphic potential $F$. 

\paragraph{Special K\"ahler domains}
Let $U \subset \bbC^{n}$ be a connected 
open domain and $F: U \rightarrow \bbC$
a holomorphic function which satisfies
the condition that the matrix
$$ \Im\left({\partial^2 F \over \partial z_i \partial z_j}\right) $$
is nondegenerate. Then the  function 
$$ k = {1 \over 2} \Im\left( \sum{\partial F \over \partial z_i} \bar{z}_i \right)\; $$
defines a K\"ahler potential on $U$.
With the corresponding K\"ahler form $ \omega = i\,  \partial \bar{\partial} k$,
and metric $g = \omega(i\, \cdot, \cdot)$, the domain $U$ is a (pseudo-) 
K\"ahler manifold\footnote{We do not require that the K\"ahler metric $g$
is definite}.
Such a domain $U$
will be called a {\em special K\"ahler domain}. On a special K\"ahler domain $U$ 
there are  flat coordinates, called {\em flat special coordinates},
\begin{equation} \label{fc}
  x_i =  \Re(z_i) \; , \; \; \;  
  y_i =  \Re({\partial F \over \partial z_i})
\end{equation}
which define 
on $U$ a torsionfree flat connection $\nabla$ so that $\omega$
is parallel. The complex manifold $U$ with the data 
$(g,\nabla )$ is then a special K\"ahler manifold. 
Conversely, 
any special K\"ahler manifold is locally equivalent 
to a special  K\"ahler domain $(U,g,\nabla )$. 

Another peculiar feature of  special K\"ahler domains
is that the K\"ahler metric $g$ is a {\em Hessian\/} metric
with respect to the flat connection.
This means that on $U$ there
exists a real potential function $f$ so that $g= \nabla df$.    
(The fact that $g$ is locally Hessian is well known. 
An explicit formula for $f$ which is given
in terms of the holomorphic function $F$, see \cite{C2}, 
shows that $f$ exists globally on $U$.)
Moreover, in the flat coordinates the smooth function $f$ satisfies the 
Monge-Amp\` ere equation
\begin{equation}  
\label{MA} 
| \det \partial^2 f |  = c \; ,  
\end{equation}
where $c>0$ is a constant.
As a consequence, the data $(g, \nabla)$ 
give $U$ the geometry of a  parabolic affine hypersphere, see
\cite{bc1}. Explicitly,  
\[\lambda(u) = (x_1(u), \ldots , x_n(u), y_1(u),  \ldots , y_n(u), f(u))\]  
defines  a Blaschke immersion $\lambda: U \rightarrow \bbR^{2n+1}$ into
affine space $\bbR^{2n+1}$ which induces the data $(\nabla, g)$.  

\subsection{The metric geometry of conic special K\"ahler 
domains} 
In this paper, we are mainly concerned with {\em conic  special K\"ahler 
domains}. We call a special  K\"ahler domain $U \subset \bbC^{n+1}\backslash \{ 0 \}$ {\em conic}, 
if $\bbC^* U \subset U$ and if the holomorphic prepotential $F$ is a 
homogeneous function of degree $2$. Moreover, we require that the
potential $k$ does not vanish on a conic special K\"ahler 
domain. Locally, any conic special K\"ahler 
manifold is equivalent to a  conic  special K\"ahler domain
$U \subset \bbC^{n+1}\backslash \{ 0 \}$.
To any conic domain $U \subset \bbC^{n+1}$ we let $\bar{U}$ 
denote its image in the projective space $\CP^n$. 
We consider the projection map 
$$  \pi:  U \longrightarrow \bar{U} $$ 
which is a submersion, and view $U$ as a principal
$\bbC^*$-bundle over $\bar{U}$. 
The special 
K\"ahler metric $g$ on $U$ naturally induces a K\"ahler metric $\bar{g}$
on $\bar{U}$ via the projection $\pi$.  
The metric $\bar{g}$ is
defined by the formula
\begin{equation} \label{pKM}
\bar{g}_{\pi(u)}(d\pi(X), d\pi(X)) = {g_u(X,X) \over g_u(u,u)}  - \left| {g_u(X,u) \over g_u(u,u)}
\right|^2 \;, \; \; X \in T_u \bbC^{n+1} \; . 
\end{equation}
(Note  
that $g$ is definite on the vertical spaces 
${\cal V}_u = \bbC u \subset T_u \bbC^{n+1}$ of the 
fibration $\pi$ by the condition that $k \neq 0$, see Lemma \ref{l1} below.) 
Let $\bar \omega$ denote the corresponding K\"ahler form on $\bar U$. Then
it is easy to see that the pull-back $\pi^* \bar{\omega}$ on $U$ is
given by 
$$ \pi^* \bar{\omega} = i \, \partial \bar{\partial} \log k \; . $$
We call the domain $\bar{U} \subset \CP^n$ with the metric $\bar{g}$ 
a {\em projective special K\"ahler domain}. 
The simplest example of such a domain is 
projective space $\CP^n$ itself with the
Fubini-Study metric: 

\begin{example} Putting $U= \bbC^{n+1} \backslash \{ 0 \}$ 
and $F(z_0, \ldots, z_n) = i \sum z_j^2$, formula (\ref{pKM}), defines
the Fubini-Study metric on $\bar{U}= \CP^n$.
The famous {\em Hopf-fibration\/} 
$$ \Sr^{2n+1} \longrightarrow \CP^n $$
exhibits the sphere $\Sr^{2n+1} = \{ u \in \bbC^{n+1} \mid |u|^2=1 \}$
as a $\Sr^1$-principal bundle over $\CP^n$. The Hopf fibration is also
known to be a Riemannian submersion with respect to the standard metric on 
the sphere if the metric on  $\CP^n$ is suitably normalized.
\end{example}
It is the content of
our next proposition that the geometric construction of the Hopf-fibration
generalizes in the context of projective special K\"ahler domains. 
To establish this result we consider now 
the K\"ahler potential $k$ on $U$. 
We remark that $k$ 
satisfies $k (\alpha u) = |\alpha|^2 k(u)$, for $\alpha \in \bbC^*$,
and, by assumption, never vanishes on $U$. We put $M_c = \{ u \in U \mid 
|k(u)|=c \}$.  
Then the level surface $M_c$ is a real hypersurface in $U \subset \bbC^{n+1}$, 
and $\Sr^1$ acts freely on $M_c$.    

\begin{proposition} \label{hyperspheres} 
The hypersurfaces $M_c \subset U$ are
nondegenerate with respect to the metric $g$.  
Moreover, $\Sr^1$ acts isometrically on $(M_c,g)$,
and $M_c$ is a $\Sr^1$-principal bundle over $\bar{U}$.
If $k>0$ then the projection map 
$$\pi_c: (M_c,g)  \longrightarrow (\, \bar{U}, \bar{g}\,)$$ 
is a semi-Riemannian submersion for $c= {1 \over 2}$. (If $k<0$ then
$\pi_c$ is an anti-isometry on horizontal vectors for $c= {1 \over 2}$)
\end{proposition} 
 
We will need a lemma. Let
$h = g - i 2 \omega$ denote the Hermitian product 
on $U$ which is defined by $g$.  We let $\xi(u)=u$ denote the position
vector field on $U$. 

\begin{lemma} \label{l1} \hspace{1cm} 
\begin{itemize} 
\item[i)]  $h(\xi,\,  \cdot \, ) =  2 \bar{\partial} k$ 
\item[ii)] $g(\xi,\,  \cdot \, ) =  dk$ 
\item[iii)]$g(\xi,\xi) = 2 k$ 
\end{itemize} 
\end{lemma} 
\begin{proof}
In the complex coordinates we have
$\xi= \sum (z_j {\partial \over \partial z_j} 
 + \bar{z}_j {\partial \over \partial \bar{z}_j})$
and 
\begin{equation} \label{h}
 h = \sum\Im\!\left({\partial^2 F \over \partial z_i \partial z_j}\right) dz_i \otimes 
d\bar{z}_j \; .
\end{equation} 
Consequently, 
\begin{equation*}
\begin{split}
    h(\xi,\, \cdot \, ) \; &= \; 
\sum \Im\!\left({\partial^2 F \over \partial z_i \partial z_j}\right) z_i d\bar{z}_j \\
     & =  \; -  {i \over 2} \left(
           \sum {\partial^2 F \over \partial z_i \partial z_j} z_i d \bar{z}_j
        -  \sum {\partial^2 \bar{F} \over \partial \bar{z}^i \partial \bar{z}^j} z_i d \bar{z}_j
                            \right) \\
     &=  \; -  {i \over 2} \left(
           \sum{\partial F \over \partial z_j} d \bar{z}_j 
        -  \sum {\partial^2 \bar{F} \over \partial \bar{z}^i \partial \bar{z}^j} z_i d \bar{z}_j
                            \right)\\
     &= \; 2 \bar{\partial} k 
\end{split}
\end{equation*}
This proves i). Now ii) follows from i) by calculating 
$$ g(\xi, \cdot) = \Re \, h(\xi, \cdot) = 
 ( \bar{\partial} k + \partial k)= dk
\; . $$ 
Equation iii) is implied by ii), taking into account that
the function $k$ is\linebreak[4] $\bbR^{>0}$-homogeneous of degree 2. 
\end{proof}

\begin{prf}{Proof of Proposition \ref{hyperspheres}}
We consider the $g$-orthogonal decomposition 
$T_u \bbC^{n+1} = {\cal V}_u \oplus {\cal H}_u$
into vertical and horizontal space which is defined by the canonical
submersion $\pi: U \rightarrow \bar{U}$. Theen 
${\cal V}_u$ is the real span of $\xi$ and $J \xi$,  
and in fact ${\cal H}_u = \{ X \in T_u \bbC^{n+1}\mid h(\xi, X) = 0 \}$.
In particular, $g(\xi, X) =0$, for $X \in {\cal H}_u$.
Therefore, by ii) from the lemma, it follows 
that ${\cal H}_u \subset \ker dk = T M_c$. 
We compute the pull back $\pi^* \bar{g}$ of the special 
K\"ahler metric $\bar{g}$ on $\bar{U}$  
on the tangent space of $M_c$. Using (\ref{pKM})
we get that $g_u(u,u) \, \pi^*\bar{g}_u  = g_u$ on $T_uM_c$.
Now, by iii) of the lemma, $g_u(u,u)= 2 k(u)$.  
The proposition follows.  
\end{prf}

\begin{proposition} \label{xi=xiProp} The vector field $\xi$, 
which is the position vector field on the conic complex domain $U$,
is also the position vector field in the affine coordinates
$x_i$, $y_i$. 
\end{proposition} 

\begin{proof} To see this,  we compute
\begin{eqnarray*}
 d {x_i}(\xi) & = & \Re \; d{z_i}(\xi) = \Re \, z_i  =x _i \\
 d {y_i}(\xi) & = & \Re \; d \! \left({ {\partial F} \over \partial z_i} \right) \! (\xi) =  
\Re \, \sum {\partial^2 F \over \partial z_i \partial z_j} z_j  = 
\Re \, {{\partial F} \over \partial z_i}  = y _i \; . 
\end{eqnarray*} 
Hence, $\xi = \sum x_i {\partial \over \partial x_i} + y_i {\partial 
\over \partial y _i}$ as claimed. 
\end{proof} 

\paragraph{Metric cones} For any manifold $M$ with a (pseudo-) Riemannian metric $g$, 
the manifold $C(M)= \bbR^{>0} \times M$ with the metric $dr^2 + r^2 g$ is
called the {\em metric cone\/} over $M$. More generally,
we consider cone metrics of the type $g_{\kappa}= {1 \over \kappa} dr^2 + r^2 g$, 
where $\kappa \neq 0$ is a constant. We denote the corresponding metric cone
as $C_{\kappa}(M)= (C(M), g_\kappa)$. Let us put ${\rm sign} \, k = 1$ if $k>0$ and
${\rm sign} \, k = -1$ if $k<0$. 

\begin{corollary} \label{metriccone} 
Let $U$ be a conic special K\"ahler domain
with K\"ahler potential $k$, and special K\"ahler metric $g$. 
Then $(U,g)$ is isometric to the metric cone 
$C_{ {\rm sign \, } k}(M_{1 \over 2})$. 
\end{corollary}
\begin{proof}
Since $\bbR^{>0}$ acts freely on $U$, the map 
$$ \Phi: C(M_{1 \over 2}) \rightarrow U \; \; (r,u) \mapsto ru$$
is a diffeomorphism. Note that $d\Phi (r{\partial \over \partial r}) = \xi$. 
The homogeneity of the holomorphic potential $F$ implies that the second derivatives
of $F$ are constant on radial lines in $U$. Hence, by formula (\ref{h}), 
we have $g_{ru}(rX, rX) =  r^2 g_u(X,X)$, for $u \in M_{1 \over 2}$, 
$X \in T_u M_{1 \over 2}$.
Moreover, by iii) of Lemma \ref{l1}, $g_{ru}(\xi, \xi)  = 2k(ru)= 
r^2 \,  {\rm sign\,  } k$.  
It is now immediate from ii) of Lemma \ref{l1} that $\Phi$ is an isometry. 
\end{proof}

\begin{proposition} The one-form $\eta := \omega(\xi, \cdot)$ 
defines a contact structure on  $M_{1 \over 2}$.
\end{proposition}
\begin{proof} $d\eta = L_\xi \omega = 2 \omega$ is nondegenerate on 
${\rm ker}\, \eta = J\xi^\perp$. 
\end{proof} 

\subsection{The affine geometry of 
conic special K\"ahler domains}  \label{localaffine}
Having just seen that any conic special K\"ahler 
domain $(U,g)$ has the geometry of a metric cone over the
level surface $(M_{1 \over 2},g)$ of $k$, we consider now the 
question
how the flat affine connection $\nabla$ on $U$ interacts with the 
cone structure of $(U,g)$. The flat affine geometry on
$U$ is determined by the coordinate change (\ref{fc}) which
embeds $U$ as a domain in $\bbR^{2n+2}$. Since the symplectic form
$\omega$ is $\nabla$-parallel, so is the
volume form $$ \theta = \theta_g \, = \; {1 \over (n+1)!} \, \; \omega^{\wedge n+1} \; .$$
Using the flat special coordinates we may view 
$$ M_{1 \over 2} \longrightarrow  \bbR^{2n+2} $$ immersed 
into affine space as a real hypersurface. In the light 
of Corollary \ref{metriccone}, the next result shows
that the metric structure on $(U,g)$ is 
determined by the affine geometry of the 
hypersurface $M_{1 \over 2}$.  

\begin{theorem} \label{Affinesphere}
In the flat special coordinates of the special K\"ahler domain
$U \subset \bbC^{n+1}$ the hypersurface $M_{1 \over 2} \subset U$ 
immerses as a non-degenerate hypersurface in $\bbR^{2n+2}$.
The transversal field $E=- {\rm sign} \, k \, \xi$ 
is a Blaschke-normal for $M_{1 \over 2}$ with respect to the volume form
$\theta$ on $\bbR^{2n+2}$, and the corresponding Blaschke-metric on $M_{1 \over 2}$
coincides with the metric $g$ induced from $(U,g)$. 
Moreover, $M_{1 \over 2}$ is an affine
hypersphere of affine mean curvature ${\rm sign} \, k$.  
\end{theorem} 

We start the proof of the theorem with a lemma. 
Any vector field $X$ on $M_{1 \over 2}$ with values in $\bbC^{n+1}$ 
has a natural extension $\tilde{X}$ on $U$ 
which is defined by $\tilde{X}(ru) = r X(u)$, for $u \in M_{1 \over 2}$.

\pagebreak[3]
\begin{lemma} \label{l2}  \hspace{1cm}
\begin{itemize} 
\item[i)]   $\xi\cdot g(\tilde{X}, \tilde{Y}) =  2 g(\tilde{X}, \tilde{Y})$,
\item[ii)]   $ (\nabla_\xi g) (\tilde{X}, \tilde{Y}) = 0$,
\item[iii)]  $ g(\xi, \nabla_{\tilde{X}} \tilde{Y}) = - g(\tilde{X}, \tilde{Y})$, if $Y$ is tangent to $M_{1 \over 2}$.
\end{itemize}
\end{lemma} 

\begin{proof}  Using Proposition \ref{xi=xiProp},  
i) follows since the function $g(\tilde{X}, \tilde{Y})$ is
$\bbR^{>0}$-homoge\-neous of degree 2. Also 
from $\tilde{X}(ru) = r \tilde{X}(u)$,
for all $u \in U$, we deduce that 
$\nabla_{\xi} \tilde{X} = \nabla_{\tilde{X}} \xi = \tilde{X}$.
Therefore, $(\nabla_\xi g) (\tilde{X}, \tilde{Y}) = \xi \cdot g (\tilde{X}, \tilde{Y})
- g (\nabla_\xi \tilde{X}, \tilde{Y})-  g (\tilde{X}, \nabla_\xi \tilde{Y})= 
 \xi \cdot g(\tilde{X}, \tilde{Y}) -2 g(\tilde{X},\tilde{Y})$. Hence, ii) follows
from i).

Now,  if $Y$ is tangent to $M_{1 \over 2}$ then 
$g(\xi, \nabla_{\tilde{X}} \tilde{Y}) + g(\tilde{X},\tilde{Y}) = 
-  (\nabla_{\tilde{X}} g) (\xi, \tilde{Y}) 
= - (\nabla_\xi g) (\tilde{X}, \tilde{Y})$, by the symmetry of $\nabla g$.
Hence, iii) follows from ii).

\end{proof}  

\begin{prf}{Proof of Theorem \ref{Affinesphere}.}
Let $X,Y$ denote vector fields tangent to $M_{1 \over 2}$, and
put $\kappa= {{\rm sign}\,  k}$. Then, by 
ii),iii) of Lemma \ref{l1}, 
and Lemma  \ref{l2} the Gau\ss-formula (\ref{Gformula}) 
for the hypersurface $M_{1 \over 2}$, 
with respect to $\xi$ reads 
$$ \nabla_X Y = \hat{\nabla}_X Y - \kappa  g(X, Y)\,  \xi \; , $$
where $\hat{\nabla}$ defines the induced connection on $M_{1 \over 2}$.
Therefore the affine metric on $M_{1 \over 2}$
with respect to the transversal vector field $E= - \kappa\,\xi$
coincides with the metric $g$. Let $\theta_{1 \over 2}$ denote the 
metric volume form of the pseudo-Riemannian manifold $(M_{1 \over 2},g)$.
To show that $E$ is a Blaschke normal, we note that (for an appropriate
choice of orientation of $M_{1 \over 2}$) the metric volume form 
$\theta = \theta_g$ of the ambient space $(U,g)$ is given by $\theta = 
 - \kappa \,dr \wedge r^{(2n+1)} \theta_{1 \over 2}$ in
the conic product coordinates $\Phi$ from the proof of Corollary \ref{metriccone}.
And, therefore, $\theta (E, \ldots ) = 
\theta( -\kappa r {\partial \over \partial r}, \ldots) = 
\theta_{{1 \over 2}}$ along $M_{1 \over 2}$. Hence, $E$ is a Blaschke-normal.
Since $A= - \nabla E = \kappa\, {\rm  Id}$, $M_{1 \over 2}$ is an affine 
hypersphere of affine mean curvature $H= \kappa$. 
\end{prf}

Now it is easy to find a $\nabla$-potential 
for $g$. 

\begin{corollary}  \label{npotential}
The  K\"ahler potential $k$ is also a $\nabla$-potential 
for the special K\"ahler metric $g$, i.e.\ $g = \nabla dk$ on $U$.  
\end{corollary} 
\begin{proof} For homogeneous vector fields $\tilde{X}$ and $\tilde{Y}$, 
we compute  
\[ (\nabla_{\tilde{X}} dk) (\tilde{Y}) = \tilde{X} \cdot d k(\tilde{Y}) - 
dk(\nabla_{\tilde{X}} \tilde{Y})\, .\]  
If $Y$ is tangent to $M_{1 \over 2}$ then,
using ii) of Lemma \ref{l1} $dk(\nabla_{\tilde{X}} \tilde{Y})=  g(\xi,\nabla_{\tilde{X}} \tilde{Y})$,
and iii) of Lemma \ref{l2} implies $(\nabla_{\tilde{X}} dk) (\tilde{Y}) = g(\tilde{X}, \tilde{Y})$.
If $\tilde{Y}= \xi$ then 
we get $(\nabla_{\tilde{X}} dk) (\xi) = \tilde{X} \cdot dk(\xi) -  dk(\tilde{X})=
\tilde{X} \cdot 2k - dk(\tilde{X}) = dk(\tilde{X}) = g(\xi, \tilde{X})$.
\end{proof}

\section{Affine Sasakian hyperspheres} \label{AS}

In Riemannian geometry a manifold $(S,g)$ is called 
Sasakian if the corresponding metric cone $(C(S), g_1)$
is a K\"ahler manifold, see e.g.\ \cite{BG}. 
More generally, we call a
(pseudo-) Riemannian manifold Sasakian if the metric 
cone $C_{\kappa}(S)$ is a \mbox{(pseudo-)} 
K\"ahler manifold. Let $U$ be a 
a conic special K\"ahler domain, and $S=M_{1\over 2} \subset U$
the affine sphere which is associated to $U$ by Theorem \ref{Affinesphere}.
By Corollary \ref{metriccone}, the affine hypersphere $S$ 
is a Sasakian manifold. 
However, the concept of Sasakian manifold does not take 
into account the presence
of the affine connection $\hat{\nabla}$ on $S$. 
Let $(S,g,\hat{\nabla})$ be a proper affine sphere.
We show 
below that the
metric cone $C_{\kappa}(S)$ admits,
as the natural affine differential geometric structure
induced from $S$, the geometry of a parabolic 
affine hypersphere 
$(C(S), g_{\kappa}, \nabla)$. This parabolic sphere is
called the {\em parabolic cone\/} over $S$. In \cite{bc1} it was remarked 
that the geometric data of a special K\"ahler manifold are
in fact the geometric data of a parabolic sphere $(M,\nabla,g)$
whose Blaschke metric is K\"ahler, and whose  K\"ahler form 
$\omega$ is $\nabla$-parallel. This motivates the following   

\begin{definition} A proper affine hypersphere $(S,g,\hat{\nabla})$ is
called an {\em affine Sasakian hypersphere\/} if the 
parabolic cone $(C(S),g,\nabla)$ over $S$ is 
K\"ahler, and the corresponding 
K\"ahler form $\omega$ is $\nabla$-parallel.  
\end{definition} 
Equivalently, a proper affine sphere $S$ is affine Sasakian, if and 
only if the parabolic cone over $S$ is special K\"ahler. 
  
\subsection{The parabolic cone over a proper affine sphere} 
We show here that every proper affine
hypersphere may be naturally realized as
a hypersurface in a conic parabolic affine sphere.
We already encountered this phenomenon,
however in the particular context of conic special K\"ahler domains. 

\paragraph{Proper spheres embed into conic parabolic spheres}
Let $(M,g)$ be a pseudo-Riemannian
manifold. We view $M = \{1\} \times M$ in a canonical 
way as a submanifold of $C_{\kappa}(M)$ with the metric $g$ 
induced from the cone metric $g_{\kappa}$ on $C(M)$. Note also 
that the multiplicative group $\bbR^{>0}$ acts 
on $C(M)$.

\begin{proposition} \label{paraboliccone} 
Let $\psi: S \longrightarrow \bbR^{n+1}$
be a proper affine hypersphere of affine mean curvature $\kappa$,
and with induced Blaschke data $(\hat{\nabla},h)$. 
Then the metric cone $C_{\kappa}(S)$
admits a torsionfree, flat, 
$\bbR^{>0}$-invariant connection $\nabla$ so that the data 
$(h_{\kappa}, \nabla)$
satisfy the integrability conditions for a
parabolic affine hypersphere. 
\end{proposition}
\begin{proof}  
We consider the local diffeomorphism 
$ \Phi: C(S) \rightarrow \bbR^{n+1}$ given by 
$(r,u) \mapsto r \psi(u)$ and let $\nabla$ 
be the pullback of the canonical flat connection on $\bbR^{n+1}$.
To simplify the notation we view $S$ as a hypersurface
in $\bbR^{n+1}$. Also we may then assume that $E= -\kappa \,\xi$ 
is the affine normal of $S$, where $\xi(x)= x$ 
is the position vector field on $\bbR^{n+1}$. 
For a vector field $X$ on $S$, let $\bar{X}$ denote
the constant extension of $X$ to the product 
manifold $C(S)= \bbR^{>0} \times S$. 
Also we define the vector field  $\tilde{X}$ on
$U= \Phi(C(S))$ by $\tilde{X}(ru) = r X(u)$, where
$u \in S$ and $r>0$. We let $\bar{\xi} = r {\partial \over \partial r}$
denote the position vector field on the cone $C(S)$.
Then $\bar{X}= \Phi^* \tilde{X}$, 
and $\bar{\xi}= \Phi^* \xi$.

We show first that the metric volume form 
$\theta_{h_\kappa}$ is $\nabla$ parallel.  
Note first that, for the right choice of orientation of $C(S)$,  
$$ \theta_{h_\kappa} \, =  \, |\kappa|^{-\frac{{1}}{2}} dr \wedge r^n \theta_h \; .$$  
We choose a (local) basis of vector fields
$X_1, \ldots, X_n$ on $S$.  
Since $\theta_h(X_1,\ldots, X_n)= \det(E,X_1, \ldots, X_n)$
along $S$, we get on $C(S)$:
$$\theta_{h_\kappa}(\bar{\xi}, \bar{X}_1, \ldots, \bar{X}_n)=
 |\kappa|^{-{1\over 2}} r^{n+1} \theta_h(X_1,\ldots, X_n) = 
 \pm |\kappa|^{1\over 2}r \det(\xi,\tilde{X}_1, \ldots, \tilde{X}_n) \; .$$ 
Therefore $\theta_{h_\kappa}= \pm |\kappa|^{1\over 2}r\Phi^* \det$, 
and hence the equiaffine condition ii) is satisfied with
respect to $\nabla$.

Next we show that $\nabla h_\kappa$ is totally symmetric. 
It is enough to verify that
$$(\nabla_X h_\kappa)(Y,Z) =(\nabla_Y h_\kappa)(X,Z) \; , $$
for all vector fields $X,Y$ and $Z$ on $C(S)$. 
We remark that if  $X,Y$ are vector fields on $S$ the following formulas 
hold on $C(S)$:
\begin{eqnarray}
\nabla_{\bar{X}}\bar{Y}   &=& \overline{\hat{\nabla}_XY}  - 
\kappa r^{-2} h_\kappa(\bar{X},\bar{Y}) \bar{\xi} \, ,  \label{wpr1}\\
\nabla_{\bar{X}}\bar{\xi} &=& \bar{X} \; , \; \;
\nabla_{\bar{\xi}}\bar{X} = \bar{X}  \label{wpr2}
\end{eqnarray}  
Therefore 
\begin{equation*} \begin{split}
(\nabla_{\bar{X}}  h_\kappa)(\bar{Y}, \bar{Z}) & = \bar{X}\cdot  h_\kappa(\bar{Y}, \bar{Z})
-  h_\kappa(\nabla_{\bar{X}} \bar{Y}, \bar{Z})-  h_\kappa (\bar{Y},\nabla_{\bar{X}} \bar{Z})\\
&= r^2 X \cdot h(Y,Z)  - r^2 h(\hat{\nabla}_X Y,Z) - r^2 h(Y, \hat{\nabla}_X Z) \; .
\end{split} \end{equation*} 
Hence, for vector fields $\bar{X},\bar{Y}, \bar{Z}$ the
compatibility condition i) for $h_\kappa$ is implied by  
i) for $h$. Next we compute
\begin{equation*} \begin{split}
(\nabla_{\bar{\xi}}  h_\kappa)(\bar{Y}, \bar{Z})&  =  
\bar{\xi} \cdot h_\kappa(\bar{Y}, \bar{Z})
-  h_\kappa(\nabla_{\bar{\xi}} \bar{Y}, \bar{Z})-  h_\kappa(\bar{Y},\nabla_{\bar{\xi}} \bar{Z}) \\
 & = 2  h_\kappa(\bar{Y}, \bar{Z}) -  h_\kappa(\bar{Y}, \bar{Z}) -  h_\kappa(\bar{Y}, \bar{Z}) =0 \; .
\end{split} \end{equation*} 
But also 
\begin{equation*} \begin{split}
(\nabla_{\bar{Y}}  h_\kappa)(\bar{\xi}, \bar{Z}) &= 
-  h_\kappa(\nabla_{\bar{Y}} \bar{\xi}, \bar{Z})-  h_\kappa(\bar{\xi},\nabla_{\bar{Y}} \bar{Z})\\
& =
-  h_\kappa({\bar{Y}}, \bar{Z}) +  
\kappa r^{-2} h_\kappa(\bar{\xi}, \bar{\xi}) h_\kappa({\bar{Y}}, \bar{Z}) =0 \;. 
\end{split} \end{equation*} 
Finally, we  easily see that 
$(\nabla_{\bar{\xi}}  h_\kappa)(\bar{X}, \bar{\xi})= 
(\nabla_{\bar{X}} h_\kappa)(\bar{\xi}, \bar{\xi})=0$.
Hence, it follows that $\nabla h_\kappa$ is totally symmetric.   
\end{proof}  

Note that $(h_\kappa, \nabla)$ satisfies the
integrability condition for {\em parabolic\/} spheres since
$\nabla$ is flat. Hence, $(C(S),h_{\kappa}, \nabla)$ 
has the structure of a parabolic affine sphere. 
As a consequence of the fundamental theorem of
affine differential geometry, 
if $C(S)$ is simply connected, the data $(h_\kappa, \nabla)$
are obtained from a Blaschke immersion 
$\Phi: C(S) \rightarrow \bbR^{n+2}$
as a parabolic affine hypersphere. 
Thus, the affine sphere $(S,h,\hat{\nabla})$ is 
realized in a canonical way as a submanifold
of a parabolic affine sphere $(C(S),h_\kappa,\nabla)$,
and the Blaschke metric on $S$, with respect to $(C(S),\nabla)$, 
coincides with the metric $h$, induced from $h_\kappa$. 
We call the parabolic affine 
sphere $(C(S),h_\kappa,\nabla)$ the {\em parabolic
cone\/} over $S$.

\paragraph{Completeness of affine spheres}
We recall an important fact about parabolic spheres. 
Calabi \cite{Cal1} and Pogorelov \cite{Pog} 
proved that if the affine metric $g$ of a parabolic affine hypersphere $(M,g,\nabla)$
is definite and complete, then $M$ must be a paraboloid. The case
that a proper affine sphere $(S,h,\hat{\nabla})$ has a definite metric 
is also of particular interest. The Blaschke normal of $S$ may be chosen so
that the affine mean curvature $H=\kappa$ 
is positive. If (with this choice of normal) the metric 
$h$ is positive definite, 
then $S$ is called an {\em elliptic\/} affine sphere, if $h$ is negative 
definite then $S$ is called {\em hyperbolic}. Therefore
$S$ is elliptic, if and only if the metric cone $C_{\kappa}(S)$ 
carries a definite metric $h_{\kappa}$. In the hyperbolic case 
the metric $h_\kappa$ has Lorentzian signature $(1,n)$. 
There is the following result of Calabi \cite{Cal2} on 
complete elliptic hyperspheres:

\begin{theorem} \label{ca2}
Let $S$ be an elliptic affine hypersphere with complete
Blaschke metric $h$. Then $S$ is an ellipsoid.
\end{theorem}

Let $S$ be an elliptic affine hypersphere with complete
metric $h$. Then the parabolic sphere 
$(C(S),h_{\kappa},\nabla)$ has definite metric $h_{\kappa}$. 
However, clearly
the metric cone $C(S)$ is not complete.
But Calabi's theorem implies that if $S$ is
complete then $C(S)= U \subset \bbR^{n+1}$ may be completed
in $0 \in \bbR^{n+1}$ to $\bar{U}= \bbR^{n+1}$, so that
the metric $h_{\kappa}$ smoothly extends to $\bbR^{n+1}$.
We deduce:

\begin{corollary} \label{ca3} 
Let $S$ be an elliptic affine hypersphere with complete
metric $h$ and affine mean curvature $\kappa$.   
Then the  parabolic cone $(C_{\kappa}(S),h_{\kappa}, \nabla)$ is obtained by  
deleting a point in an elliptic paraboloid. 
\end{corollary}

\subsection{Characterization of affine Sasakian hyperspheres} \label{ASH} 
Let $(S,g)$ be a (pseudo-) Riemannian manifold, $D$ 
the Levi-Civita connection on $S$. Then 
a {\em Sasakian structure\/} on $S$ is provided by a Killing 
vector field $\sigma$ of constant length 
$g(\sigma , \sigma )= \kappa^{-1}$ so that the covariant
derivative $\Phi = D \sigma$ satisfies 
\[    (D_X\Phi) (Y) = \kappa g(\sigma ,Y) X - g(X,Y) \sigma \; .\]
The Killing vector field $\sigma$ and the one-form
$\eta  = \kappa \, g(\sigma , \cdot)$ are called the 
{\em characteristic vector field\/} and the
{\em characteristic one-form\/} of the Sasakian
structure on $S$. Let $C_{\kappa}(S)$ be a metric
cone over $S$, and let $\xi = r {\partial \over \partial r}$
denote the Euler field  on $C(S)$. We define 
a complex structure $J$ on $C(S)$ by the formulas
$$  J \bar{X} = \overline{\Phi X} - \eta(X) \xi \; , \; J \xi = \sigma \; .$$
It is straightforward to verify that in fact $J^2 = -{\rm Id}$, 
and that the cone metric $g_{\kappa}$ 
is $J$-invariant. Moreover $J$ is parallel with respect to the Levi-Civita
connection.  Hence, $J$ is integrable and $C_{\kappa}(S)$ is K\"ahler. 
Conversely, if $C_{\kappa}(S)$ is K\"ahler with
respect to the complex structure $J$ then $\sigma = J  \xi$ 
defines the characteristic vector field of a Sasakian
structure on $S$.

\begin{proposition} \label{ASc}
Let $(S,g,\hat{\nabla})$ be a proper 
affine hypersphere with Sasakian structure $\sigma $. 
Then the parabolic cone over $(S,g,\hat{\nabla})$ is 
special K\"ahler with respect to the complex structure
$J$ induced from $\sigma $ if and only if $\Phi = \hat{\nabla} \sigma $. 
\end{proposition} 

\begin{proof} Let us first recall the formulas (\ref{wpr1}), 
(\ref{wpr2}) from the proof of Proposition 
\ref{paraboliccone}, which are satisfied by the flat connection 
$\nabla$ on
$C(S)$. Note also that the same (warped product) 
relations hold for the metric connections $D$ and $\bar{D}$,
where $\bar{D}$ is the Levi-Civita connection of the
cone metric $g_\kappa$. Next we remark that the parabolic cone 
$(C(S),g_\kappa,\nabla)$
is special K\"ahler if and only if the special K\"ahler condition
\begin{equation} \label{special}
d^{\nabla} J = 0
\end{equation} 
is satisfied. For a vector field $Y$ on $S$, we
compute $(\nabla_{\xi} J) \bar{Y} =0$, and 
$(\nabla_{\bar{Y}} J) \xi = \nabla_{\bar{Y}} \bar{\sigma } - J \bar{Y}$, where
$J \bar{Y} = \overline{D_Y \sigma } - \eta(Y) \xi$ and 
$\nabla_{\bar{Y}} \bar{\sigma } = \overline{\hat{\nabla}_Y \sigma } -  \kappa g(Y,\sigma ) \xi$. 
Therefore if (\ref{special}) is satisfied $(\nabla_{\bar{Y}} J) \xi =0 $,
and hence $\hat{\nabla}_\cdot \sigma  = D_\cdot \sigma =\Phi$.
Conversely, from $\Phi= \hat{\nabla}_\cdot \sigma$ we deduce that 
$J \bar{Y} = \nabla_{\bar{Y}} \bar{\sigma}$
and hence, since $\nabla$ is flat, it follows  (\ref{special})
along $S$. Moreover, from the above equations 
$d^{\nabla}\! J(\bar{Y}, \xi) = (\nabla_{\bar{Y}} J) \xi=0$ follows
immediately. Therefore,
the parabolic cone  $(C(S),g_\kappa,\nabla)$ is special K\"ahler.
\end{proof}  

Consequently, if the Sasakian structure $\sigma$ satisfies 
$\Phi = \hat{\nabla} \sigma$ we call $\sigma$ an {\em affine Sasakian
structure\/} on the hypersphere $(S,g,\hat{\nabla})$.


\section{Applications} \label{Applications}

\subsection{The Canonical circle bundle}
Let $\pi: M \rightarrow \bar{M}$ be a projective special K\"ahler manifold,
where the conic manifold $M$ carries the data $(J,g,\nabla)$.
Let $(\tilde{M},J,g,\nabla)$ be the universal covering
space of $M$, and $\lambda: \tilde{M} \rightarrow V$ a compatible 
Lagrangian embedding into a pseudo-Hermitian, 
symplectic vector space $(V,\gamma, \Omega)$. 
Since the embedding $\lambda$ is unique up to isometry 
of  $(V,\gamma, \Omega)$, the function 
$$    k(p) \, = \, {1 \over 2} \, \gamma(\, \lambda(p),\, \lambda(p)) \; ,\, \;
p \in \tilde{M}  $$  
is invariant under deck-transformations of the covering, and hence defines 
a function $k: M \rightarrow \bbR^{>0}$. Note that, by iii) of Lemma \ref{l1} 
and by Corollary \ref{npotential}, $(M,\nabla,g)$ 
is a Hessian-manifold with potential $k$. We define a 
family of hypersurfaces $M_c= \{ p\in M \mid  k(p)  = c\}$
in $M$. Then the hypersurfaces $M_c$ are invariant 
by the natural isometric $S^1 \subset \bbC^*$ action on
the conic manifold $M$. We call $$ S := M_{\frac{1}{2}} \longrightarrow \bar{M}$$
the {\em canonical circle bundle} over the projective special K\"ahler manifold $\bar{M}$. 
\begin{theorem} \label{mainThm} 
Let $\bar{M}$ be a projective special K\"ahler manifold
and $S \rightarrow \bar{M}$ its canonical
circle bundle. Then $S$ has a canonical structure of 
a proper affine hypersphere. Moreover, $S$ carries 
an affine Sasakian structure which determines the
projective special K\"ahler geometry on $\bar{M}$.  
\end{theorem}
\begin{proof} It is enough to prove the theorem locally. 
Therefore we assume $M_{1 \over 2} \subset U$, where
$U$ is a special K\"ahler domain with data 
$(g,J,\nabla)$. By Theorem \ref{Affinesphere}, $S=M_{1 \over 2} \subset U$
is a proper affine sphere, so that $(U,g)$ is the metric cone 
over $S$. Since the flat coordinates on $U$ are conic, i.e.\  
$\bbR^{>0}$-equivariant, the flat connection $\nabla$ on $U=C(S)$ coincides with
the flat connection on $C(S)$ which is constructed in Proposition \ref{paraboliccone}.
Hence, $(U,g,\nabla)$ is the parabolic cone over $S$, and the parabolic
cone is special K\"ahler.  
In particular, the sphere $S$ is affine Sasakian, and, by Proposition \ref{ASc},
the Sasakian structure $\sigma$ on $S$ induced from $J$ is affine Sasakian. 
 \end{proof} 


\subsection{Projective special K\"ahler domains with a definite metric}
Let $\bar U$ be a projective special K\"ahler domain with
a definite metric $\bar{g}$ and $F$ the potential function
of the corresponding  special 
K\"ahler domain $U \subset \bbC^{n+1} \backslash \{0\}$
which carries the special K\"ahler metric $g$ defined by formula (\ref{h}).
Note that by formula (\ref{pKM}) the function $-F$ 
induces the same metric $\bar{g}$ on $\bar{U}$, however 
the signature of the metric $g$ on $U$ is inverted. 

\begin{definition} 
A projective special K\"ahler domain $\bar{U}$ with
a definite metric $\bar{g}$ is called of
{\em elliptic type\/} if the metric $g$ on $U$ 
is definite.
\end{definition}  

We remark that if $\bar{U}$ is an elliptic 
projective special K\"ahler domain, then by 
formula (\ref{pKM})  the metric $\bar{g}$ on 
$\bar{U}$ must be positive definite. Moreover 
the affine hypersphere $S \subset U$ which
is associated to $\bar{U}$ by Theorem 
\ref{Affinesphere} has a definite metric, and 
$S$ is an elliptic affine hypersphere. 
Conversely, if $\bar{U}$ is a projective special K\"ahler domain
with a negative definite metric $\bar{g}$,  then
the associated affine hypersphere $S$
has an affine metric with Lorentzian 
signature.

\paragraph{Characterization of complex projective space} \label{CPn}
In \cite{Lu} it was proved that a special K\"ahler 
manifold $M$ with a (positive) definite complete metric
is flat. In fact, it may also be deduced from this 
result that any complete special K\"ahler domain 
$U \subset \bbC^{n}$ with a definite metric is just $\bbC^{n}$ with
a Hermitian inner product. In the case of  
projective special K\"ahler domains there 
are many (homogeneous) examples with a definite and complete 
metric known, for instance, the examples 
given in the section \ref{homex}. Among elliptic special K\"ahler 
domains though,  the projective space $\CP^n$ is 
characterized by its completeness property:

\begin{theorem} Let $\bar{U} \subset \CP^n$ be a
projective special K\"ahler domain of elliptic
type with a complete metric $\bar{g}$. Then 
$\bar{U}=\CP^n$ and $\bar{g}$ is homothetic
to the Fubini-Study metric on $\CP^n$. 
\end{theorem}
\begin{proof} 
We may choose $F$ on $U \subset \bbC^{n+1}$ 
so that $g$ is positive definite. Therefore
the  K\"ahler potential $k$ on $U$ is positive. By 
Theorem \ref{Affinesphere}, the associated
affine hypersphere $S$ is of elliptic type
with a positive definite metric and, 
since $S \longrightarrow \bar{U}$ is a 
Riemannian submersion with a complete 
base and compact fibre $S^1$, $S$ has a  
complete metric as well. Hence $S$ is an ellipsoid by
Thm \ref{ca2}. Recall that, by Corollary
\ref{metriccone}, the special K\"ahler domain 
$U \subset \bbC^{n+1}$ over $\bar{U}$ 
is the parabolic cone over $S$ and, 
by Corollary \ref{ca3}, $U= \bbC^{n+1} \backslash \{0\}$.
Also by  Corollary \ref{ca3}, the metric $g$ on 
$U$ has a quadratic potential with respect to the
flat connection $\nabla$ on $U$. Since, by Corollary \ref{npotential},
$k$ is a $\nabla$-potential for $g$, $k$ must
be a homogeneous quadratic function in the affine 
coordinates. Hence, it follows that 
the cone metric $g$ is parallel 
with respect to $\nabla$. Therefore $\nabla=D$, 
which is possible only if $F$ is a quadratic
function and $g$ is just a Hermitian inner 
product on $\bbC^{n+1}$. In this case, $\bar{g}$
is homothetic to the Fubini-Study metric. 
\end{proof} 

\subsection{Calabi-Yau moduli space} \label{CalabiYau}
We recall that a {\em Calabi-Yau $m$-fold} (of general type) is an 
oriented compact Riemannian manifold $(X,g)$ with holonomy group 
${\rm Hol}(X,g) = {\rm SU}(m)$. This implies that $X$ admits a unique
complex structure $J$ compatible with the orientation such that
$(X,J,g)$ is a K\"ahler manifold and a parallel 
$J$-holomorphic $(m,0)$-form ${\rm vol}$ (a holomorphic volume form),  
which is unique up to constant scale. 
In particular, $(X,J)$ is a complex manifold of (complex)  
dimension $m$ with trivial canonical bundle $\wedge^{m,0}T^*X$. Let
$\bar{M}$ be the Kuranishi moduli space of $(X,J)$, i.e.\ the (local) 
moduli space of complex structures $I$ on $X$. There is a natural 
holomorphic line bundle over $\bar{M}$ whose fibre at $I \in \bar{M}$ is  
$\Gamma_{hol} (\wedge_I^{m,0}T^*X) = H^{m,0}(X,I)$ ($\Gamma_{hol}$ stands
for {\em holomorphic} sections).   Let $\pi : M \rightarrow 
\bar{M}$ be the corresponding holomorphic\linebreak[3] $\bbC^*$-bundle:  
$\pi^{-1}(I) = H^{m,0}(X,I) - \{ 0\}$. The one-dimensional complex vector 
spaces  $H^{m,0}(X,I)$ have a natural norm: 
$\| {\rm vol} \|^2 := (\sqrt{-1})^{-m} \int_X 
{\rm vol}\wedge \overline{\rm vol}$.
Let $S \subset M$ be the unit circle bundle with respect to that norm.  
     
\begin{theorem} 
Let $S \rightarrow \bar{M}$ be the above circle bundle over the Kuranishi
moduli space of a Calabi-Yau threefold. Then $S$ has naturally the
structure of a Lorentzian affine Sasakian hypersphere. In particular,
$S$ is a proper affine hypersphere.   
\end{theorem} 
\begin{proof} It is known that $\bar{M}$ has the structure of a projective 
special K\"ahler manifold. We briefly recall the construction of that
structure. (For more details, see \cite{C1}). 
The cup product defines a  complex symplectic form $\Omega$ on 
$V := H^3(X,\bbC)$ and  $\gamma = \sqrt{-1}\Omega (\cdot , \bar{\cdot})$ is
a pseudo-Hermitian form of (complex) signature $(n+1,n+1)$, where 
$n = h^{1,2} = \dim \bar{M}$. The map 
\[ \bar{M} \ni I \mapsto H^{3,0}(X,I) \in P(V) \] 
is a holomorphic immersion and is 
induced by a conic holomorphic 
immersion $\phi : M \rightarrow V -\{ 0\}$, with 
the following properties: $\phi^*\Omega = 0$ ($\phi$ is Lagrangian) and
$g = {\rm Re} \, \phi^*\gamma$ is a K\"ahler metric of complex signature $(1,n)$ 
on the complex manifold $M$. These properties correspond to the
first and second Hodge-Riemann bilinear relations for the underlying
variation of Hodge structure of weight $3$. As explained in section \ref{SK} 
the conic immersion $\phi$ induces on $M$ the structure of a conic 
special K\"ahler manifold such that the corresponding 
projective special K\"ahler metric on $\bar{M}$ is negative definite 
(according to the conventions of this paper). Moreover, the circle
bundle $S$ defined above coincides with the canonical circle bundle
$S = M_{\frac{1}{2}}$ of the projective special K\"ahler manifold $\bar{M}$
(notice that $(\sqrt{-1})^{-m} = \sqrt{-1}$ for $m=3$ and hence
$\|u\|^2 = \gamma (u,u)$ for $u\in H^{3,0}(X,I)$).      
Now we can apply Theorem \ref{mainThm}. 
\end{proof}


\section{Homogeneous examples} \label{homex}
The basic example of an affine Sasakian hypersphere $S$ is provided
by the total space of the Hopf fibration 
\[ S = S^{2n+1} = {\rm SU}(n+1)/{\rm SU}(n)\longrightarrow \CP^n = 
{\rm SU}(n+1)/{\rm S (U}(n){\rm U}(1))\; .\] 
In the Lagrangian picture the corresponding conic affine special K\"ahler manifold 
$(M,J,g,\nabla )$ is given 
as a linear Lagrangian subspace 
$M \subset V=T^*\bbC^{n+1}$ for which 
the restriction of the Hermitian metric $\gamma$ is positive definite.   
Since $M$ is a linear subspace the flat connection $\nabla$ 
coincides with the Levi-Civita connection $D$ of
$g = \Re \, \gamma$.  
The group ${\rm SU}(n+1)$ acts transitively
on $\bar{M} = \CP^n$ by holomorphic isometries of the special K\"ahler metric 
(Fubini-Study metric). The action is induced from 
the canonical linear symplectic action of ${\rm SU}(n+1)$
on $V=T^*\bbC^{n+1}$ which preserves the Hermitian metric $\gamma$ and the Lagrangian
subspace $M \subset V$. This action preserves also
the  affine Sasakian hypersphere $S^{2n+1} \subset M$ and induces a 
transitive action on $S^{2n+1}$ preserving the affine geometric 
and Sasakian structures. 

More generally, one can consider Lagrangian
subspaces $M \subset V=T^*\bbC^{n+1}$ of arbitrary Hermitian signature
$(p,q)$, $p+q = n+1$.  They correspond
to fibrations 
\[ S = {\rm SU}(p,q)/{\rm SU}(p,q-1)\longrightarrow 
{\rm SU}(p,q)/{\rm S (U}(p,q-1){\rm U}(1)) = \bar{M}\, .\] 
The case $q=1$ is of particular interest. In that case 
the projective special K\"ahler metric is negative definite 
(as for the Calabi-Yau moduli space and as for the target manifolds 
of N=2 D=4 supergravity theories with vector multiplets) 
and hence the metric of the affine Sasakian hypersphere has 
Lorentzian signature: 
$\bar{M} = \CH^n$ is the complex hyperbolic space and $S$ 
is the real hyperbolic (2n+1)-space of Lorentzian signature 
(anti de Sitter space).

\paragraph{The Classification}
A  projective special K\"ahler manifold $\bar{M}=P(M)$ will be called
{\em homogeneous\/} if it admits a transitive group of isometries $G$
whose action is induced by a $G$-action on the conic manifold $M$ 
preserving the data $(g,J,\nabla)$. 
Homogeneous projective special K\"ahler manifolds 
$$ \bar{M} = P(M) = G/K$$ 
with $K$ compact have been classified in \cite{AC} 
under the assumption that $G$ is a real semisimple Lie group. 
We recall the result only in the most interesting 
case of negative definite metric on $\bar{M}$. It turns out
that in this case 
the manifolds $\bar{M}= G/K$ are Hermitian symmetric 
spaces of non-compact type
and are in one-to-one correspondence with the complex simple Lie algebras
${\fr l}$ different from ${\fr c}_n = {\fr sp}(\bbC^{2n})$. 
In all the cases the underlying
conic affine special K\"ahler manifold is a  
Lagrangian cone $M \subset V$ 
generated by the $G$-orbit of a 
highest weight vector of a $G^{\bbC}$-module $V$ of symplectic type. 
The $G^{\bbC}$-module $V$ admits a $G$-invariant 
real structure $\tau$ compatible
with the symplectic structure $\Omega$, which defines a
Hermitian metric  $\gamma = \sqrt{-1}\Omega (\cdot , \tau \cdot )$. 
The affine special
K\"ahler metric is the restriction of $g = \Re \, \gamma$ to $M$. 
The list is the 
following:
\begin{enumerate}
\item[A)] ${\fr l} ={\fr {sl}}_{n+3}(\bbC ), \quad 
\bar{M} = \CH^{n} = {\rm SU}(n,1)/{\rm S(U}(n){\rm U}(1)), \quad 
V = \bbC^{n+1} \oplus (\bbC^{n+1})^*$  
\item[BD)] ${\fr l} = {\fr {so}}_{n+5}(\bbC ), \quad 
\bar{M} = ({\rm SL}(2,\bbR )/{\rm SO}(2)) \times 
({\rm SO}(n-1,2)/{\rm SO}(n-1){\rm SO}(2)),\linebreak[4] V = 
\bbC^2\otimes \bbC^{n+1}$ 
\item[E6)] ${\fr l} ={\fr e}_6(\bbC ), \quad 
\bar{M}  = {\rm SU}(3,3)/{\rm S(U}(3) {\rm U}(3)), \quad 
V = \bigwedge^3\bbC^6$
\item[E7)] ${\fr l} ={\fr e}_7(\bbC ), \quad
\bar{M}  = {\rm SO}^*(12)/{\rm U}(6), \quad V^{(32)} = V(\pi_6)$ (semispinor) 
\item[E8)] ${\fr l} ={\fr e}_8(\bbC ), \quad
\bar{M}  = {\rm E}_7^{(-25)}/{\rm E}_6 {\rm SO}(2), \quad V^{(56)} 
= V(\pi_1)$
\item[F)] ${\fr l} ={\fr f}_4(\bbC ), \quad
\bar{M}  = {\rm Sp}(\bbR^6)/{\rm U}(3), \quad V^{(14)}(\pi_3) = 
\bigwedge_0^3\bbC^6$ 
\item[G)] ${\fr l} ={\fr g}_2(\bbC), \quad 
\bar{M}  = \CH^1 = {\rm SL}(2,\bbR )/{\rm SO}(2), \quad 
V = \bigvee^3\bbC^2.$
\end{enumerate}

\noindent 
Here $V(\lambda )$ denotes the irreducible module $G^{\bbC}$-module with
highest weight $\lambda = \sum \lambda_i \pi_i$, where $\pi_i$ are 
the fundamental weights. The notation  $V^{(d)}$ indicates that the
module has complex dimension $d$.   
Notice that in the cases A) and BD) $n = \dim_{\bbC} \bar{M}$. 
The only redundancy in this list occurs the case $n=1$. In fact, the
Dynkin  diagrams $A_3 = B_3$ define the same projective special
K\"ahler manifold $\CH^1 = {\rm SU}(1,1)/{\rm S(U}(1){\rm U}(1)) 
= {\rm SL}(2,\bbR )/{\rm SO}(2)$. In both cases the corresponding conic manifold 
$M$ is a linear Lagrangian subspace in the vector space $V$.  

Note that it
may happen that projective special K\"ahler manifolds are  
isometric as Riemannian manifolds, but nevertheless their special
geometry is different: The diagram $G_2$ defines 
$\bar{M}  = \CH^1$ but in this case  the underlying
conic affine special K\"ahler manifolds  $M  \subset V$ is not
a linear subspace, as for type A), n=1. 
In fact, $V = \bigvee^3\bbC^2$ is the symmetric cube
of the defining representation $\bbC^2$ of $G^{\bbC} = {\rm SL}(2,\bbC )$. 
The Zariski closure of $M \subset V$ is the nonlinear cone 
$M' = \{ u^3| u\in \bbC^2 \} \subset V$ and $M \subset M'$ is open. 

\paragraph{Homogeneous affine Sasakian spheres}
An affine hypersphere $(S,g,\hat{\nabla} )$ is called {\em homogeneous} if 
${\rm Aut}(S) = {\rm Aut}(S,g,\hat{\nabla} )$ acts transitively on $S$.  
Note that in general ${\rm Aut}(S)$ is a proper subgroup of
${\rm Isom}(S) = {\rm Aut}(S,g)$. If $S$ has
an affine Sasakian structure $\sigma$ then let ${\rm Aut}_{\sigma}(S)$
be the subgroup of those automorphisms in ${\rm Aut}(S)$ which commute with 
the flow of the vector field $\sigma$. We call ${\rm Aut}_{\sigma}(S)$ 
the group of automorphisms of the affine Sasakian sphere $S$.
Clearly, any affine Sasakian hypersphere $S$ with a transitive 
action of ${\rm Aut}_{\sigma}(S)$ is a circle bundle over a 
homogeneous projective special K\"ahler manifold. If $\bar{M}$ is
homogeneous then ${\rm Isom}(S)$ acts transitively on $S$. 
But note that, in 
general, the canonical isometric $S^1$-action on $S$ does not
preserve the connection $\hat{\nabla}$.
The following theorem is a 
consequence of the above classification. 
\begin{theorem} \label{homo1} Let $S$ be the affine Sasakian
hypersphere over a homogeneous projective special K\"ahler manifold 
$\bar{M} = G/K$ of a real semisimple Lie group $G$. 
If the special K\"ahler metric of $\bar{M}$ is negative
definite then $\bar{M}$ belongs to the above list A)-G) and 
$G$ acts transitively by automorphisms 
of the Lorentzian affine Sasakian hypersphere $S$.   
\end{theorem} 

\begin{proof} By construction, the $G$-action on
$\bar{M}$ is induced by a $G$-action on the symplectic
vector space $V$ which preserves the geometric data on $V$.
Hence $G$ acts also on the canonical circle bundle $S$
over $M$ preserving the affine Sasakian geometry on $S$.
Note now that in all the cases the centre $Z(K) \cong {\rm U}(1)$ of 
$K$ acts non-trivially, and hence transitively, on 
the fibre of $S\rightarrow \bar{M}$ over the canonical base point
$o = eK$ in $\bar{M} = G/K$. This follows, for example, from the fact that
$K$ is the stabilizer of the line $l = {\bbC}v \subset V$ generated by a 
highest weight vector $v\in V^{\tau}$ of the $G^{\bbC}$-module $V$. In fact, 
$K$ contains a (compact) Cartan subgroup of $G$, which cannot act trivially
on $l$. (Notice, that the semisimple part of $K$, however,  
acts trivially on $l$.)    
\end{proof}

Clearly, the ${\rm Aut}(S)$ action on the affine sphere 
$S$ extends to a linear (with respect to the flat connection $\nabla$)
action on the parabolic cone which contains $S$. 
Hence, if  $G \subset {\rm Aut}(S)$ acts transitively, the 
affine sphere $S$ arises as a generic $G$-orbit in 
a real vector space $W$. If $G$ is
semisimple then $S$ must be the level-set of a 
homogeneous $G$-invariant polynomial on $W$. 

\begin{theorem} \label{homo2}
Let $S$ be an affine Sasakian hypersphere with 
Lorentzian metric. If ${\rm Aut}_{\sigma}(S)$ contains 
a semisimple  transitive group $G$ then the affine
sphere $(S,g,\hat{\nabla})$ arises  
as a hypersurface which
is defined by a $G$-invariant homogeneous quartic polynomial
on a real vector space $W$. 
\end{theorem} 

\begin{proof} $S$ identifies with the canonical 
circle bundle in the parabolic cone $M=C(S)$ which is 
special K\"ahler. The action of  ${\rm Aut}_{\sigma}(S)$ on $S$ 
extends to an action on $M$ which preserves the special 
K\"ahler data on $M$. Using a compatible Lagrangian immersion
we may therefore as well assume 
that the action of $G={\rm Aut}_{\sigma}(S)$ on $S$ is 
induced by an action of $G$ on a Hermitian symplectic 
vector space $(V,\gamma,\Omega)$. 
In fact, we identify $S$ as an affine sphere in the real 
vector space $W=V^{\tau}$, and $S$ is a level set of the K\"ahler potential 
$k$, which is, as a function on $V^{\tau}$, homogeneous of degree 2
and invariant by $G$. 
We claim
that $k^2$ is a quartic polynomial. Since $G$ acts with cohomogeneity one,
it is sufficient to show that $V^\tau$ admits a homogeneous
$G$-invariant quartic polynomial, which is then necessarily proportional
to $k^2$. To show this it is clearly enough to construct 
a (complex) homogeneous $G^{\bbC}$-invariant quartic polynomial on $V$. 
The existence of such a polynomial on $V$ follows  by the 
following general argument.  

As we know, the $G^{\bbC}$-module $V$ is associated to
a Dynkin diagramm $\Delta$ of the type A, B, D, E, F or G.
We give some more detail how this correspondence works. 
(See \cite{AC} for a complete account.) 
Let $N = N(\Delta ) = L/L_o$ be the compact symmetric 
quaternionic K\"ahler manifold which is associated to the 
Dynkin diagramm $\Delta$. (See \cite{Wolf}.) 
$L$ is the compact simple Lie group with trivial
centre associated to $\Delta$ and 
$L_o = {\rm Sp}(1)H$ is the
stabilizer of a point $o\in N$. The
complexified isotropy representation is  
a product $T_o N\otimes {\bbC}  = \bbC^2\otimes_{\bbC} V$. 
The group ${\rm Sp}(1)$ acts by the
standard representation on $\bbC^2$, and $V$ is a complex module for $H$ 
which admits a skew symmetric bilinear invariant. 
It follows that the maximal semisimple 
subgroup $H' \subset H$ is a compact form of a complex semisimple
group $G^{\bbC}$ which acts on $V$. In this way, we have associated a 
$G^{\bbC}$-module $V$ to the Dynkin diagramm $\Delta$. 
Now the quaternionic Weyl tensor, see \cite{Sa}, of the 
quaternionic K\"ahler manifold $N = N(\Delta )$ at the  point $o\in N$ 
gives rise to a nonzero $G^{\bbC}$-invariant 
element $Q \in S^4V^*$. 
This shows the existence of a nontrivial
$G^{\bbC}$-invariant homogeneous quartic polynomial $Q$ on $V$.  
\end{proof} 

In examples it is not difficult to guess the 
quartic invariant $Q$  directly from the $G$-module $V^{\tau}$.
This gives an explicit description of the corresponding 
affine hyperspheres.\\

\noindent {\bf Examples}\\

\noindent 
A) $G = {\rm SU}(n,1)$,  $V^{\tau} = \bbC^{n,1}$, $Q(v) = g(v,v)^2$, 
where $g$ is the ${\rm SU}(n,1)$-invariant Hermitian product.\\

\noindent BD) $G = {\rm SL}(2,\bbR ) \times {\rm SO}(n-1,2))$, $V^{\tau} = 
\bbR^2 \otimes \bbR^{n-1,2} \cong {\rm Hom}((\bbR^2)^{*},\bbR^{n-1,2})$. 
Let $\omega$ be a
${\rm SL}(2,\bbR )$-invariant symplectic form on $\bbR^2$ and
$g$ the ${\rm SO}(n-1,2)$-invariant scalar product, defining 
identifications $\Phi_\omega: \bbR^2 \cong (\bbR^2)^*$, 
$\Phi_g: \bbR^{n-1,2} \cong (\bbR^{n-1,2})^*$. For $A \in 
{\rm Hom}((\bbR^2)^{*},\bbR^{n-1,2})$ let 
$A^* \in  {\rm Hom}((\bbR^{n-1,2})^*,\bbR^2)$
be the dual morphism. Then 
$Q(A) = \det (A^* \Phi_g   A \Phi_\omega )$.\\  
 
\noindent E6) $G^{\bbC} = {\rm SL}(6,\bbC)$, $V = (\wedge^3\bbC^6)^*$. 
To any 3-form $\alpha$ we can
associate the operator 
\[ A_{\alpha} :  \bbC^6 \, \longrightarrow\,    
\left(\bigwedge^5 \bbC^6\right)^* = \, \bbC^6\, ,\quad 
v \mapsto \alpha \wedge \iota_v\alpha\, .\]
Then $Q(\alpha ) = {\rm trace} (A_{\alpha}^2)$. It is easy to check that
$Q\neq 0$ by evaluating $Q$ on $dz^1\wedge dz^2\wedge dz^3 + dz^4\wedge
dz^5\wedge dz^6$. This example is discussed in detail in \cite{H} and  
the corresponding real symplectic ${\rm SL}(6,\bbR )$-module is also
considered. Here we are interested in the real structure $\tau$ invariant
under the real form $G = {\rm SU}(3,3)$ of ${\rm SL}(6,\bbC)$. It is 
induced by the ${\rm SU}(3,3)$-invariant pseudo-Hermitian form  
on $\bbC^6 = \bbC^{3,3}$. In fact, this form induces a $G$-invariant 
pseudo-Hermitian form $\gamma$ on $V = (\wedge^3\bbC^6)^*$. This 
determines a $G$-invariant real structure $\tau$ on $V$ such that
$-i\gamma (\cdot ,\tau \cdot ) = \Omega$ is the $G^{\bbC}$-invariant
symplectic form of $V$: $\Omega (\alpha ,\beta )
dz^1\wedge dz^2\wedge dz^3\wedge dz^4 \wedge dz^5\wedge dz^6 = 
\alpha \wedge \beta$. \\ 
  
\noindent 
F)  $G = {\rm Sp}(\bbR^6)$,  $V^{\tau} = \bigwedge_0^3\bbR^6$ is the 
kernel of the map $\bigwedge^3\bbR^6 \ni \alpha \mapsto \omega \wedge \alpha
\in \bigwedge^5\bbR^6$, where $\omega$ is the symplectic form on $\bbR^6$. 
The $G$-invariant quartic polynomial $Q$ is just the restriction of the 
${\rm SL}(6,\bbC )$-invariant quartic polynomial on $\bigwedge^3\bbC^6$, 
see previous example, to the subspace $\bigwedge_0^3\bbR^6$. \\

\noindent    
G) $G = {\rm SL}(2,\bbR )$, $V^{\tau} = \bigvee^3\bbR^2$. The elements
of $V^{\tau}$ can be considered as homogeneous cubic polynomials $p$ on
$\bbR^2$. Let $q(p)= \det (\partial^2 p) \in \bigvee^2\bbR^2$ be the
determinant of the Hessian of $p$. Then $Q = D(q(p))$ is the discriminant of 
$q(p)$.
 
\medskip
\noindent 
Remarks:\\   
1) In all the above examples (A-E) the
group $R^*\cdot G$ acts with an open orbit on $V^{\tau}$, in other 
words $V^{\tau}$ with the
action of $R^*\cdot G$ is a real {\it prehomogeneous vector space}. Complex
irreducible prehomogeneous vector spaces were classified in \cite{SK}.\\
2) The Sasaki field $\sigma$ of the affine hypersphere $S \subset V^{\tau}$ 
can be easily computed from the real quartic invariant $Q$. {}From Lemma
\ref{l1} ii) it follows that $\sigma = J\xi$ is precisely 
the Hamilton vector field $X_k$ associated to the K\"ahler potential $k$. 
We can normalize the $G$-invariant real symplectic structure on $V^{\tau}$ 
(or the invariant $Q$) such that $Q$ is related to the K\"ahler potential $k$
by the formula $Q = k^2$. Then we have $X_Q = 2kX_k$ and therefore, 
since $k = 1/2$ on $S$, we have $\sigma = X_Q$ on $S$.

\paragraph{Compact quotients} 
Let $G$ be one of the real semi-simple Lie groups from the 
list A)-G). By Theorem \ref{homo1}, $G$ acts transitively and
properly on a Lorentzian affine hypersphere $S= G/\tilde{K}$
which fibers over a Hermitian symmetric space $G/K$ of non-compact
type, $K = \tilde{K} {Z}(K)$.   
By a result of Borel \cite{Borel}, $G$ 
admits cocompact lattices $\Gamma \leq G$. 
This allows to construct compact Clifford-Klein 
forms $$     S_{\Gamma}  =  \; {_{\displaystyle \Gamma}}    \lmod \, G/ {{\displaystyle \tilde{K}}}  $$
for the Lorentzian homogeneous spaces $S= G/\tilde{K}$. 
The spaces $S_{\Gamma}$ admit an isometric $S^1$-action 
(induced from the affine Sasakian structure) with finite stabilizers, the orbit space being a 
Hermitian locally symmetric space $$ \bar{M}_{\Gamma}
= \; {_{\displaystyle \Gamma}}    \lmod \,  G/K \; .$$
In his influential paper \cite{Kul}, Kulkarni observed the 
existence of non-trivial circle bundles over compact locally complex 
hyperbolic spaces, carrying a Lorentzian metric of constant curvature 1.
This corresponds to the  complex hyperbolic case 
$\bar{M} = \CH^n = \SU(n,1)/S(\U(n) \U(1))$, i.e.\ case A) in our list. 
In this sense, our construction generalizes Kul\-kar\-ni's 
construction of  {\em compact\/} Lorentzian space-forms. 
It seems worthwile to further study the 
particular Lorentzian geometry of the homogeneous spaces $S$ 
occuring in examples B) to G), and 
their compact Clifford-Klein forms.    
However, in this paper we content ourselves with summarizing what
was just explained:   

\begin{corollary} Let $\bar{M}$ be one of the Hermitian symmetric
spaces appearing in the list A)-G), and $\bar{M}_{\Gamma}$ a compact
Clifford-Klein form for $\bar{M}$. Then $\bar{M}_{\Gamma}$ is the orbit
space of an isometric $S^1$-action on a compact Clifford-Klein $S_{\Gamma}$
for the Lorentzian homogeneous space $S$ associated to $\bar{M}$. 
\end{corollary}

\end{document}